\documentclass[a4paper,12pt]{article}

\usepackage{amssymb}
\usepackage{amsmath}
\usepackage{mathptmx}
\usepackage{color}
\usepackage{cancel}
\usepackage{bbm}
\usepackage{pstricks}
\usepackage{psfrag}
\def\ap{'\thinspace}

\def\spacingset#1{\def\baselinestretch{#1}\small\normalsize}
\spacingset{1.3}

\newtheorem{lemma}{Lemma}[section]
\newtheorem{theorem}{Theorem}[section]
\newtheorem{remark}{Remark}[section]

\newtheorem{proposition}{Proposition}[section]

\def\be{\begin{equation}}
\def\ee{\end{equation}}
\def\bea{\begin{eqnarray}}
\def\eea{\end{eqnarray}}
\def\beann{\begin{eqnarray*}}
\def\eeann{\end{eqnarray*}}
\def\bsea{\begin{subeqnarray}}
\def\esea{\end{subeqnarray}}
\def\bmat{\left[ \begin{array}}
\def\emat{\end{array} \right]} 

\def\ap{'\thinspace}

\def\proof{\noindent{\bf{\em Proof:}\ \ }}
\def\QED{\mbox{\rule[0pt]{1.5ex}{1.5ex}}}
\def\endproof{\hspace*{\fill}~\QED\par\endtrivlist\unskip}
\newcommand{\real}{{\mathbb{R}}}

\def\gR{{\cal R}}
\def\gS{{\cal S}}

\def\gV{{\cal V}}

\newcommand{\ima}{\operatorname{im}}

\newcommand{\defi}{\stackrel{\text{\tiny def}}{=}}
\def\tra{{\scalebox{.6}{\mbox{T}}}}
\def\bsmat{\left[ \begin{smallmatrix}}
\def\esmat{\end{smallmatrix} \right]}

\definecolor{Royalblue}{cmyk}{1,0.30,0.2,0.2}

\begin{document}
\begin{titlepage}
\title{\vspace{-10mm}
The generalised continuous algebraic Riccati equation \\
and impulse-free continuous-time LQ optimal control\thanks{Partially supported by the Italian Ministry for Education and Research (MIUR) under PRIN grant n. 20085FFJ2Z ``New Algorithms and Applications of System Identification and Adaptive Control" and by the Australian Research Council under the grant FT120100604. Research carried out while the first author was visiting Curtin University, Perth (WA), Australia. }\vspace{10mm}}
\author{{\large Augusto Ferrante$^\ddagger$  \quad Lorenzo Ntogramatzidis$^\star$ }\\
       {\small $^\ddagger$Dipartimento di Ingegneria dell\ap Informazione}\\[-3mm]
        {\small    Universit\`a di Padova, via Gradenigo, 6/B -- 35131 Padova, Italy}\\[-3mm]
       {\small     {\tt augusto@dei.unipd.it}} \\ 
       {\small     $^\star$Department of Mathematics and Statistics}\\[-3mm]
       {\small     Curtin University, Perth WA, Australia.}\\[-3mm]
       {\small    {\tt L.Ntogramatzidis@curtin.edu.au}}
 }%
\thispagestyle{empty} \maketitle \thispagestyle{empty}

\begin{center}
\begin{minipage}{14.8cm}
\begin{center}
\bf Abstract
\end{center}
The purpose of this paper is to investigate the role that the continuous-time generalised Riccati equation plays within the context of singular linear-quadratic optimal control. This equation has been defined following the analogy with the discrete-time generalised Riccati equation, but, differently from the discrete case, to date the importance of this equation in the context of optimal control is yet to be understood. This note addresses this point. 
We show in particular that when  the continuous-time generalised Riccati equation admits a symmetric solution, the corresponding linear-quadratic (LQ) problem admits an impulse-free optimal control.
\end{minipage}
\end{center}
\begin{center}
\begin{minipage}{14.2cm}
\vspace{2mm}
\end{minipage}
\end{center}
\thispagestyle{empty}
\end{titlepage}
%
\section{Introduction}
\label{secintro} 
Riccati equations are universally regarded as a
cornerstone of modern control theory. 
In particular, it is well known that the solution of the classic finite and infinite-horizon LQ optimal control problem strongly depends on the matrix weighting the input in the cost function, traditionally denoted by $R$. When $R$ is positive
definite, the problem is said to be {\it regular}
(see e.g. \cite{Anderson-M-89,Lewis-S-95}), whereas when $R$ is positive semidefinite, the problem is called
{\it singular}. The singular cases have been treated
within the framework of geometric control theory, see for example
\cite{Hautus-S-83,Willems-KS-86,Saberi-S-87} and the references cited therein. In particular, in
\cite{Hautus-S-83} and \cite{Willems-KS-86} it was proved that an optimal solution of the singular LQ problem exists for all initial conditions if the class of allowable controls is extended to include
distributions. \\[-0.6cm]

In the discrete time,  the solution of finite and infinite-horizon LQ problems can be found resorting to the so-called  {\em generalised discrete algebraic Riccati equation}. In particular, the link between the solutions of LQ problems and the solutions of generalised discrete algebraic/difference equations have been investigated in \cite{Rappaport-S-71,Ferrante-N-12-sub1} for the finite horizon and in \cite{Ferrante-N-12-sub} for the infinite horizon.  A similar generalisation has been carried out for the continuous-time algebraic Riccati equation in \cite{Ionescu-O-96-1}, where the {\em 
generalised Riccati equation} was defined in such a way that the inverse of $R$ appearing in the standard Riccati equation is replaced by its pseudo-inverse. Some conditions under which this equation admits a stabilising  solution were  investigated  in  terms 
of  the so-called deflating subspaces of the extended  Hamiltonian  pencil. 
Some preliminary work on the continuous-time algebraic Riccati equation within the context of spectral factorisation has been carried out in \cite{Chen-F-89} and \cite{Weiss-94}. Nevertheless, to date the role of this equation in relation to the solution of optimal control problems in the continuous time {has not been fully} explained. 
The goal of this paper is to fill this gap, by providing a counterpart of the results in  \cite{Ferrante-N-12-sub} for the continuous case. In particular, we describe the role that the generalised continuous algebraic Riccati equation plays in singular LQ optimal control.
Such role does not trivially follow from the analogy with the discrete case. Indeed, in the continuous time, whenever the optimal control involves distributions, none of the solutions of the generalised Riccati equation is optimising. The goal of this paper is to address this delicate issue. 
Thus, the first aim of this paper is to explain the connection of the generalised continuous-time algebraic Riccati equation and of the 
generalised Riccati differential equation -- which is also defined by substitution of the inverse of $R$ with the pseudo-inverse -- and the 
 solution of the standard LQ optimal
control problem with infinite and finite horizons, respectively.
We will show that when the generalised Riccati equation possesses a symmetric solution, both the finite and the infinite-horizon LQ problems admit an impulse-free solution. Moreover, such control can always be expressed as a state-feedback, where the gain can be obtained from the solution of the generalised continuous-time algebraic/differential Riccati equation. 
   We also provide an insightful geometric characterisation of the situations in which the singular LQ problem admits impulse-free solutions, in terms of the largest output-nulling controllability subspace and of the smallest input-containing subspace of the underlying system.\\[-5mm]

{\bf Notation.}  The image and the kernel of matrix $M$ are denoted by $\ima\,M$ and $\ker\,M$, respectively, while the transpose and the Moore-Penrose pseudo-inverse of $M$ are denoted by $M^\tra$ and $M^\dagger$, respectively.

 \section{Generalised CARE}
{ We consider} the following matrix equation
\begin{equation}
 X\,A+A^\tra\,X-(S+X\,B)\,R^{\dagger}\,(S^\tra\!+B^\tra X)+Q=0, \label{gcare}
\end{equation}
with $Q, A\in \real^{n \times n}$, $B,S \in \real^{n \times m}$, $R \in \real^{m \times m}$ and we make the following standing assumption:
\begin{equation}
\label{sd}
\Pi \defi \bmat{cc}  Q & S \\ S^\tra & R \emat=\Pi^\tra \ge 0.
\end{equation}
Eq. (\ref{gcare}) is often referred to as the {\em generalised continuous algebraic Riccati equation} GCARE($\Sigma$), and represents a generalisation of the classic continuous algebraic Riccati equation arising in infinite-horizon LQ problems since here $R$ is allowed to be singular. Eq. (\ref{gcare}) along with the additional condition
\begin{equation}
 \ker R \subseteq \ker (S+X\,B), \label{kercond}
\end{equation}
 will be referred to as {\em constrained generalised continuous algebraic Riccati equation}, and is denoted by CGCARE($\Sigma$). Observe that from (\ref{sd}) we have $\ker R \subseteq \ker S$, which implies that (\ref{kercond}) is equivalent to $\ker R \subseteq \ker (X\,B)$. 
 
 The following notation is used throughout the paper. First, let $G \defi I_m-R^\dagger R$ be  the orthogonal projector that projects onto $\ker R$. Hence, $R^\dagger R$ is the orthogonal projector that projects onto $\ima R^\dagger=\ima R$. In fact, $\ker R=\ima G$. Moreover, we consider a non-singular matrix $T=[T_1\mid T_2]$ where
$\ima T_1=\ima R$ and $\ima T_2=\ima G$, and we define $B_1\defi BT_1$ and $B_2 \defi BT_2$. Finally, to any $X=X^\tra \in \real^{n \times n}$  we associate\begin{eqnarray}
Q_X& \defi & Q+A^\tra X+X\,A, \qquad
S_X   \defi   S+X\, B, \label{defgx} \\
K_X & \defi & R^\dagger\, (S^\tra+B^\tra\,X)=R^\dagger\, S_X^\tra, \qquad A_X \defi  A-B\,K_X, \\
  \Pi_X  & \defi & \left[ \begin{array}{cc} Q_X & S_X \\ S_X^\tra & R \end{array} \right]. \label{KX}
 \end{eqnarray}

\begin{lemma}
\label{lemm1}
Let $X$ be a solution of CGCARE($\Sigma$). Then, $X\,B_2=0$.
\end{lemma}
\proof From (\ref{kercond}) and from $\ker R=\ima G$, we get $(S+X\,B)\,G=0$. Moreover, since $\Pi\geq 0$, $\ker S \supseteq \ker R$. This means that $K\in \real^{n \times m}$ exists such that $S=K\,R$. Therefore, 
$S\,R^\dagger\,R=K\,R\,R^\dagger=K\,R=S$, so that $S\,G=S-S\,R^\dagger\,R=0$. Hence, $\ima (B\,G)\subseteq \ker X$ or, equivalently, $XB_2=0$.
\endproof
 
\begin{lemma}
\label{lem0}
Let $F=A-B\,R^\dagger S^\tra$ and $\Lambda=Q-S\,R^\dagger S^\tra$. Then, $\Lambda \ge 0$ and GCARE($\Sigma$) defined in (\ref{gcare}) has the same set of symmetric solutions of the following equation:
\begin{equation}
\label{gcare1}
X\,F+F^\tra\,X-X\,B\,R^{\dagger}\,B^\tra X+\Lambda=0
\end{equation}
\end{lemma}
\proof
Matrix $\Lambda$ is the generalised Schur complement of $R$ in $\Pi$. Therefore, since $\Pi$ is positive semidefinite, such is also $\Lambda$. 
The rest of the proof is a matter of standard substitutions of $F$ and $\Lambda$ into (\ref{gcare1}) to verify that (\ref{gcare}) is obtained. 
\endproof

\begin{remark}
{\em The result established for GCARE($\Sigma$) in Lemma \ref{lem0} extends without further difficulties to the  so-called generalised Riccati differential equation GRDE($\Sigma$) 
 \begin{eqnarray}
&& \dot{P}(t)+P(t)\,A+A^\tra\,P(t) \nonumber \\
&& \qquad \qquad -(S+P(t)\,B)\,R^{\dagger}\,(S^\tra\!+B^\tra P(t))+Q=0. \label{grde} 
\end{eqnarray}
 Indeed, we easily see that (\ref{grde}) has the same set of symmetric solutions of the equation:
\begin{equation}
\dot{P}(t)+P(t)\,F+F^\tra\,P(t)-P(t)\,B\,R^{\dagger}\,B^\tra P(t)+\Lambda=0. 
\end{equation}
}
\end{remark}

\begin{lemma}
Let $X=X^\tra$ be a solution of CGCARE($\Sigma$). 
Let $\gR(F,B_2)$ be the reachable subspace of the pair $(F,B_2)$. Then
\begin{eqnarray*}
&&{\bf (1)} \;\, \ker X \subseteq \ker \Lambda; \\
&& {\bf (2)} \;\, X\,\gR(F,B_2)=\{0\}; \\
&&  {\bf (3)}\;\, \Lambda\,\gR(F,B_2)=\{0\}.
\end{eqnarray*}
\end{lemma}
\proof
{\bf (1).} Let $\xi \in \ker X$. Multiplying (\ref{gcare1}) to the left by $\xi^\tra$ and to the right by $\xi$, we get $\xi^\tra\,\Lambda\,\xi=0$. 
{Since  $\Lambda\geq 0$, this implies $\Lambda\,\xi=0$. Hence,} $\ker X \subseteq \ker \Lambda$. \\
{\bf (2).} Let $\xi \in \ker X$. Post-multiplying (\ref{gcare1}) by $\xi$ we find $X\,F\,\xi=0$. This implies that $\ker X$ is $F$-invariant.
In view of Lemma \ref{lem0}, the subspace $\ker X$ contains $\ima B_2$. Hence, it contains $\gR(F,B_2)$ that is the 
smallest $F$-invariant subspace containing $\ima B_2$.
%
  This implies $\gR(F,B_2) \subseteq \ker X$. \\
  {\bf (3).} This follows directly from the chain of inclusions $\gR(F,B_2) \subseteq \ker X \subseteq \ker \Lambda$.
  \endproof


 \section{The finite-horizon LQ problem}

 \begin{lemma}
 \label{leml}
 Let $H = H^\tra\ge0$ be such that $H\, \gR(F,B_2)=\{0\}$. 
   If CGCARE($\Sigma$) (\ref{gcare}-\ref{kercond}) admits solutions,
  the generalised Riccati differential equation 
  \begin{eqnarray}
&& \dot{P}_T(t)+P_T(t)\,A+A^\tra\,P_T(t) \nonumber \\
&& \qquad \qquad -(S+P_T(t)\,B)\,R^{\dagger}\,(S^\tra\!+B^\tra P_T(t))+Q=0, \label{grde1} 
\end{eqnarray}
 with the terminal condition 
 \begin{equation}
 \label{number}
 P_T(T)=H
 \end{equation}
admits a unique solution for all $t \le T$, and this solution satisfies $P_T(t)\,B\,G=0$ for all $t \le T$.
 \end{lemma}
 \proof
Consider a set of coordinates in the input space such that the first coordinates span $\ima R$ and the second set of coordinates spans $\ima G=\ker R$. In this basis $R$ can be written as $R=\left[ \begin{smallmatrix} R_1 & 0\\[1mm] 0&0\end{smallmatrix} \right]$ with $R_1$ being invertible. In the same basis, matrix $B$ can be partitioned accordingly as $B=[\,B_1\;\;B_2\,]$ as shown above, so that $\ima B_2=\ima (B\,G)$.
Let us now consider the change of basis matrix $U=[\,U_1\;\;U_2\,]$ in the state space where $\ima U_1=\gR(F,B_2)$, so that
\begin{eqnarray*}
U^{-1} F\,U&=&\bmat{cc} \!\! F_{11} \!\!  & F_{12} \!\!  \\  \!\! O  \!\! & F_{22}  \!\! \emat \!\! , \quad U^{-1} B_1=\bmat{cc}  \!\! B_{11} \!\!  \\  \!\! B_{12}  \!\! \emat \!\! , \quad
 U^{-1} B_2=\bmat{cc}  \!\! B_{21} \!\!  \\ \!\!   O  \!\! \emat,
\end{eqnarray*}
and $ U^\tra \Lambda U=\bsmat O & O \\[1mm] O & \Lambda_{22} \esmat$
where we have used the fact that $\Lambda\,\gR(F,B_2)=\{0\}$. In this basis, since we are assuming $H\, \gR(F,B_2)=\{0\}$, we can write $U^\tra\,H\,U=\bsmat O & O \\[1mm] O & H_{22} \esmat$. 
 Consider the following matrix function $P_T(t)=\bsmat O & O \\[1mm] O & P_{22}(t) \esmat$, 
 where $P_{22}(t)$ satisfies
\begin{eqnarray}
 && \dot{P}_{22}(t)\!+\!P_{22}(t) F_{22}\!+\!F_{22}^\tra P_{22}(t)\!-\!P_{22}(t) V P_{22}(t)\!+\!\Lambda_{22}\!=\!0 \label{grde122} \\
&& P_{22}(T)=H_{22}, \label{condin22}
\end{eqnarray}
 in which $V$ is the sub-block 22 of the matrix $B\,R^\dagger B^\tra$. From \cite[Corollary 2.4]{Freiling-JS-00} we find that, since $\Pi=\Pi^\tra \ge 0$ and $H = H^\tra\ge0$, both (\ref{grde1}) and (\ref{grde122}) admit a unique solution defined in $(-\infty,T]$. 
 It is easy to see that $P_T(t)=\bsmat O & O \\[1mm] O & P_{22}(t) \esmat$, where $P_{22}(t)$, $t \in (-\infty,T]$, is the solution of (\ref{grde122}-\ref{condin22}), solves (\ref{grde1}) and (\ref{number}).
 We can therefore conclude that $P_T(t)$ is the unique solution of (\ref{grde1}-\ref{number}).
  Moreover, this solution satisfies $P_T(t)\,B_2=0$ for all $t \le T$ since in the chosen basis
$P_T(t)\,B_2=\bsmat O & O \\[1mm] O & P_{22}(t) \esmat \bsmat  B_{21} \\[1mm]  O \esmat=0$.
  \endproof
 \ \\[-0.9cm]
 
 Now we consider the generalised Riccati problem GRDE($\Sigma$) (\ref{grde1}-\ref{number})
 in relation with the finite-horizon LQ problem, which consists in the minimisation of the performance index
 \begin{eqnarray}
 J_{T,H}(x_0,u)&=&\int_0^T \bmat{cc} x^\tra(t) & u^\tra(t) \emat \bmat{cc} Q & S \\ S^\tra & R \emat \bmat{c} x(t) \\ u(t) \emat\,dt, \nonumber \\
 && \quad \qquad +x^\tra(T)\,H\,x(T)  \label{cost1}
 \end{eqnarray}
 where we only assume $\Pi=\bsmat Q & S \\[1mm] S^\tra & R \esmat=\Pi^\tra\ge 0$ subject to
 \begin{equation}
 \label{eqsys}
 \dot{x}(t)=A\,x(t)+B\,u(t), 
 \end{equation}
  and the constraint on the initial state $x(0)=x_0 \in \real^n$.   The following theorem is the first main result of this paper. It shows that when CGCARE($\Sigma$) admit a solution, the finite-horizon LQ problem always admits an impulse-free solution.
 
 \begin{theorem}
 \label{the31}
 Let CGCARE($\Sigma$) admit a solution. 
 The finite-horizon LQ problem (\ref{cost1}-\ref{eqsys}) admits impulse-free optimal solutions. All such solutions are given by
 \begin{equation}
 \label{optcontr}
 u(t)=-(S^\tra\!+\!B^\tra P_T(t))\,R^\dagger \,x(t)+G\,v(t),
 \end{equation}
 where $v(t)$ is an arbitrary regular function, and $P_T(t)$ is the solution of (\ref{grde1}) with the terminal condition (\ref{number}). The optimal cost is $x_0^\tra\,P_T(0)\,x_0$.
 \end{theorem}
 \proof 
  Let us first assume that $H\, \gR(F,B_2)=\{0\}$. The cost (\ref{cost1}) can be written for any matrix-valued differentiable function $P(t)$ as 
\begin{eqnarray*}
J_{T,H}(x_0,u)&=&\int_0^T \bmat{cc} x^\tra(t) & u^\tra(t) \emat \bmat{cc} Q & S \\ S^\tra & R \emat \bmat{c} x(t) \\ u(t) \emat\,dt \\
&& +x^\tra(T)\,H\,x(T)+\int_0^T\dfrac{d}{dt} \left( x^\tra(t)\,P(t)\,x(t)\right)\,dt\\
&& +x^\tra(0)\,P(0)\,x(0)-x^\tra(T)\,P(T)\,x(T) \\
&& \hspace{-2cm}= \!\! \int_0^T \!\!\! \bmat{cc} \!\! x^\tra(t) & u^\tra(t)  \!\! \emat  \!\!\!\! \bmat{cc} \!\!  Q\!+\!A^\tra P(t)\!+\!P(t) A\!+\!\dot{P}(t) \!& P(t) B\!+\!S\!\! \\ S^\tra\!+\!B^\tra P(t) & R \emat  \!\! \!\!\! \bmat{c} \!\!  x(t) \!\!  \\  \!\! u(t) \!\!  \emat \! dt\\
&& +x^\tra(T)\,(H-P(T))\,x(T)+x^\tra(0)\,P(0)\,x(0),
\end{eqnarray*}
where we have used the fact that
\begin{eqnarray*}
\frac{d}{dt} \left( x^\tra(t)\,P(t)\,x(t)\right) &=& \dot{x}^\tra(t)\,P(t)\,x(t) \\
&& +x^\tra(t)\,\dot{P}(t)\,x(t)+x^\tra(t)\,P(t)\,\dot{x}(t).
\end{eqnarray*} 
  Let us now consider $P(t)=P_T(t)$ to be the solution of (\ref{grde1}) with final condition $P_T(T)=H$. Since, as proved in Lemma \ref{leml},  the identity $P_T(t)\,B\,G=0$ holds for all $t \le T$, we have also $\ker (P_T(t)\,B)\subseteq \ker R$ for all $t \le T$. Thus, $\ker (P_T(t)\,B+S) \subseteq \ker R$ for all $t \le T$, and we can write
 \begin{eqnarray*}
&& \bmat{cc}  Q+A^\tra P_T(t)+P_T(t) A+\dot{P_T}(t)   &   P_T(t) B+S   \\   S^\tra+B^\tra P_T(t)   &   R   \emat \\
&& \quad =
  \bmat{cc}  (P_T(t) B+S)\,R^\dagger\,(S^\tra+B^\tra P_T(t))    &   P_T(t) B+S   \\   S^\tra+B^\tra P_T(t)   &   R   \emat
 \\ 
&& \quad =  \bmat{cc}   O   &   (P_T(t) B+S) R^\dagger R^{\frac{1}{2}}   \\   O   &   R^{\frac{1}{2}}   \emat
 \bmat{cc}   O\phantom{R^{\frac{1}{2}}}   &   O   \\  R^{\frac{1}{2}} R^\dagger (S^\tra+B^\tra  P_T(t))    &   R^{\frac{1}{2}}   \emat
 \end{eqnarray*}
 since $\ker (P_T(t)\,B+S) \subseteq \ker R$ gives $(P_T(t) B\!+\!S)\,R^\dagger R=(P_T(t) B\!+\!S)$.  Hence,
 \begin{eqnarray*}
 J_{T,H}(x_0,u) &=& \int_0^T || R^{\frac{1}{2}} R^\dagger (S^\tra+B^\tra P_T(t)) \,x(t) + R^{\frac{1}{2}}\,u(t) ||_2^2 \,dt\\
 && +x^\tra(0)\,P_T(0)\,x(0),
\end{eqnarray*}
since $P_T(T)=H$. If there exists a control $u(t)$ for which 
\begin{equation}
\label{set}
R^{\frac{1}{2}} R^\dagger (S^\tra+B^\tra P_T(t)) \,x(t) + R^{\frac{1}{2}}\,u(t)=0
\end{equation}
for all $t \in [0,T)$, then the cost function is minimal in correspondence with this control and all minimising controls satisfy (\ref{set}). The set of controls satisfying (\ref{set}) can be parameterised as $u(t)=-R^\dagger (S^\tra+B^\tra P_T(t)) \,x(t)+G\,v(t)$, where $G=I_m-R^\dagger R$ and $v(t)$ is arbitrary.\footnote{Note that this has to be understood in a $L_2$ sense.}
Now, consider the case $\gR(F,B_2) \nsubseteq \ker H$. Consider the change of basis $U=[\,U_1\;\; U_2\,]$ where $\ima U_1=\gR(F,B_2)$, and where $\ima U_2$ is the orthogonal complement of $\ima U_1$. Changing coordinates gives
$\bsmat U_1^\tra \\[1mm] U_2^\tra \esmat H \bsmat U_1 & U_2 \esmat=\bsmat H_{11} & H_{12} \\[1mm] H_{12}^\tra & H_{22} \esmat$.
 Let us now perform a further change of coordinates with the matrix $\bsmat I & U_{21} \\[1mm] O & I \esmat$ such that $U_{21}=-H_{11}^\dagger\,H_{12}$. There holds 
 \begin{equation*}
 \bmat{cc} I & O \\[0mm] U_{21}^\tra & I \emat \bmat{cc} H_{11} & H_{12} \\[1mm] H_{12}^\tra & H_{22} \emat \bmat{cc} I & U_{21}\\[1mm] O & I \emat=
 \bmat{cc} H_{11} & O \\[1mm] O & \tilde{H}_{22} \emat,
 \end{equation*}
  where $\tilde{H}_{22} \defi H_{12}^\tra\,U_{21}+H_{22}$. Thus, by writing the cost function in this new basis we get
 \begin{eqnarray*}
 J_{T,H}(x_0,u) &=&\int_0^T \bmat{cc} \!\!  x^\tra(t) \!\!& u^\tra(t)\!\! \emat\!\! \bmat{cc}\!\! Q \!\!& S\!\! \\ \!\!S^\tra\!\! & R\!\! \emat\!\! \bmat{c}\!\! x(t)\!\! \\ u(t) \!\! \emat\,dt \\
 &&+
 \bmat{cc}\!\! x_1^\tra(T)\!\! &  x_2^\tra(T)\!\! \emat\!\! \bmat{cc}\!\! H_{11} \!\!& O\!\! \\ \!\! O\!\! & \tilde{H}_{22} \!\!\emat\!\!\bmat{c}\!\! x_1(T)\!\! \\  x_2(T)\!\! \emat.
 \end{eqnarray*}
 Clearly, $J_{T,H}(x_0,u) \ge J_{T,H^\prime}(x_0,u)$, 
 where $H^\prime= \bsmat O & O \\[1mm] O & \tilde{H}_{22} \esmat$. On the other hand, in the optimal control defined by $J_{T,H^\prime}(x_0,u)$ there is a degree of freedom which is the component of the state trajectory on $\gR(F,B_2)$. In other words, in the minimisation of $J_{T,H^\prime}(x_0,u)$ we can decide to drive $x_1$ to the origin in $[0,T]$ without destroying optimality. As such, $J_{T,H}(x_0,u) = J_{T,H^\prime}(x_0,u)$. The penalty matrix of the final state in this new performance index satisfies 
$\bsmat O & O \\[1mm] O & \tilde{H}_{22} \esmat \gR(F,B_2)=\{0\}$.
 \endproof

\begin{remark}
{\em Notice that when $H=0$, we have $P_T(0) \le P_{T+\delta T}(0)$ for all $\delta T\ge 0$, because $J^\ast_{t,0}(x_0,u)$ is a non-decreasing function in $t$.}
\end{remark}

We are now interested in studying $P_T(0)$ when the terminal condition vanishes, i.e., when $H=0$, and the time interval increases. To this end, we consider a generalised Riccati differential equation where the time is reversed, and where the terminal condition becomes an initial condition, which is now equal to zero. More specifically, we consider the new matrix function 
$X(t)=P_t(0)=P_T(T-t)$.
We re-write GRDE($\Sigma$)  as a differential equation to be solved forward:
\begin{eqnarray}
&& \dot{X}(t) =X(t)\,A+A^\tra\,X(t) \nonumber \\
&& \qquad \qquad -(S+X(t)\,B)\,R^{\dagger}\,(S^\tra+B^\tra X(t))+Q, \label{grde111} \\
&& X(0)=0. \label{condin111}
\end{eqnarray}

 In the following theorem, the second main result of this paper is introduced. This theorem determines when the infinite-horizon LQ problem admits an impulse-free solution, and the set of optimal controls minimising the infinite-horizon cost
 \begin{eqnarray}
 \label{costinf}
J_\infty(x_0,u)&=&\int_0^\infty \bmat{cc} x^\tra(t) & u^\tra(t) \emat \bmat{cc} Q & S \\ S^\tra & R \emat \bmat{c} x(t) \\ u(t) \emat\,dt
\end{eqnarray}
subject to the constraint (\ref{eqsys}).

 \begin{theorem}
 Suppose CGCARE($\Sigma$) admits symmetric solutions, and that for every $x_0$ there exists an input $u(t) \in \real^m$, with $t\ge0$, such that
$J_\infty(x_0,u)$ in (\ref{costinf}) is finite. Then we have:
\begin{description} 
 \item{\bf (1)} A solution $\bar{X}=\bar{X}^\tra\ge 0$ of CGCARE($\Sigma$) is obtained as the limit of the time varying matrix generated by integrating (\ref{grde111}) with the zero initial condition (\ref{condin111}).
 \item{\bf (2)} The value of the optimal cost is {\em $x_0^\tra \bar{X} x_0$}.
 \item{\bf (3)} $\bar{X}$ is the minimum positive semidefinite solution  of CGCARE($\Sigma$).
 \item{\bf (4)} The set of {\em all} optimal controls minimising $J_\infty$ in (\ref{costinf}) can be parameterised as 
 \be
u(t)=-R^\dagger S_{\bar{X}}^\tra\,x(t)+G\,v(t), \label{optcontrinf}
\ee
with arbitrary $v(t)$.
\end{description}
\end{theorem}
\proof {\bf (1).} Consider the problem of minimising 
\begin{equation*}
J_{t,0}=\displaystyle\int_{0}^{t} [\, x^\tra(\tau) \;\;\; u^\tra(\tau) \,] \Pi \bmat{c} x(\tau) \\[0mm] u(\tau) \emat d\tau
\end{equation*}
 subject to (\ref{eqsys}) with assigned initial state $x_0 \in \real^n$.
 From Theorem \ref{the31} the optimal control for this problem exists, and the optimal cost is equal to $J_\infty^\ast(x_0,u)=x_0^\tra\, P_T(0)\,x_0=x_0^\tra\, X(T)\,x_0$. We have already observed that $X(t)=P_t(0)$ is an increasing flow of matrices in the sense of the positive semidefiniteness of symmetric matrices, i.e., $X(t+\delta t) \ge X(t)$ for all $\delta t\ge 0$.  We now show that $X(t)$ is bounded.  Indeed, given the $i$-th canonical basis vector $e_i$ of $\real^n$, we have $e_i^\tra X(t) \,e_i\le J_{\infty}(e_i,\bar{u}_i)$, where $\bar{u}_i$ is a control that renders $J_{\infty}(e_i,\bar{u}_i)$ finite (which exists by assumption). Thus,
 \[
 0 \le X(t) \le I_n\cdot \max\left\{ J_{\infty}(e_i,\bar{u}_i):\, i \in \{1,\ldots,n\} \right\} \quad \forall \, t\ge 0.
 \]
 Therefore, $X(t)$ is bounded. Taking the limit on both sides of (\ref{grde111}) we immediately see that $\bar{X} \ge 0$ is indeed a solution of CGCARE($\Sigma$).
 
{\bf (2).}  Let \\[-5.0mm]
\be
\label{infimu}
J^\circ(x_0) \defi \inf_{u} J_{\infty}(x_0,u).
\ee
Clearly, $J^\circ(x_0)\geq J_t^*(x_0)=x_0^\tra X(t)\, x_0$ for all $t \ge 0$. Then, by taking the limit, we get $J^\circ(x_0)\geq x_0^\tra \bar{X} x_0$.
We now show that the time-invariant feedback control $u^*_t \defi -K_{\bar{X}} x_t$, where $K_{\bar{X}}=R^\dagger\, (S^\tra+B^\tra\,\bar{X})$,
 yields the cost $x_0^\tra \bar{X} x_0$, which is therefore the optimal value of the cost.
Consider the cost index $J_{T,\bar{X}}(x_0,u)$.
The optimal cost for this index is achieved by using the controls satisfying (\ref{optcontr}), where $P_T(t)$ is constant and equal to $\bar{X}$, since $\bar{X}$ is a stationary solution of (\ref{grde1}-\ref{number}) and $H=\bar{X}$. Therefore, an optimal control for this index is given by the time-invariant feedback $u^*(t)=-K_{\bar{X}} x(t)$. The optimal cost does not depend on the length $T$ of the time interval and is given by 
$J_{T,\bar{X}}^*=x_0^\tra \bar{X} x_0$. 
Now we have
\begin{eqnarray}
 x_0^\tra \bar{X} x_0 & \leq & J^\circ(x_0)\leq J(x_0,u^*) \nonumber \\
 &=&\int_{0}^\infty \bmat{cc} x^\tra(t) & (u^*)^\tra(t) \emat \Pi \bmat{c} x(t) \\ u^*(t) \emat\,dt    \nonumber   \\
&=&\lim_{T\rightarrow\infty} 
\int_{0}^T \bmat{cc}  x^\tra(t)   & (u^*)^\tra(t)    \emat \Pi \bmat{c}    x(t)    \\    u^*(t)    \emat dt \nonumber \\
&=&
\lim_{T\rightarrow\infty} J_{T,\bar{X}}^* -x_T^\tra \bar{X} x_T\leq
\lim_{T\rightarrow\infty} x_0^\tra \bar{X} x_0 =x_0^\tra \bar{X} x_0. \label{chainofine-eq}
\end{eqnarray}
Comparing the first and last term of the latter expression we see that all the inequalities are indeed equalities, so that the infimum in (\ref{infimu}) is  a minimum and its value is indeed $x_0^\tra \bar{X} x_0$.

{\bf (3).} 
Suppose by contradiction that there exist another positive semidefinite solution $\tilde{X}$ of CGCARE($\Sigma$) and a vector $x_0\in\real^n$ such that
$x_0^\tra\,\tilde{X}\,x_0 <  x_0^\tra\,\bar{X}\,x_0$. 
Take the time-invariant feedback $\tilde{u}(t)=-K_{\tilde{X}} x(t)$.
The same argument that led to (\ref{chainofine-eq}) now gives 
$J_{\infty}(x_0,\tilde{u})\leq x_0^\tra \tilde{X} x_0<  x_0^\tra\,\bar{X}\,x_0$, which is a contradiction because we have shown that $  x_0^\tra\,\bar{X}\,x_0$ is the optimal value of $J_{\infty}(x_0,u)$.
 
{\bf (4).}  Consider 
 \begin{eqnarray*}
 J_\infty^\ast &=& x_0^\tra \,\bar{X} \,x_0 \\
 &=& \inf_{u(t),\, t \ge 0} \left[ \int_0^\infty \bmat{cc} x^\tra(t) & u^\tra(t) \emat \bmat{cc} Q & S \\ S^\tra & R \emat \bmat{c} x(t) \\ u(t) \emat\,dt \right].
 \end{eqnarray*}
 This infimum is indeed a minimum because we know that the optimal cost $x_0^\tra \,\bar{X} \,x_0$ can be obtained for some $u^\ast$. Hence, $J_\infty^\ast = \min_{u(t),\, t \ge 0} \left[ \displaystyle\int_0^\infty [\, x^\tra(t) \;\;\; u^\tra(t) \,] \bsmat Q & S \\[1mm] S^\tra & R \esmat \bsmat x(t) \\[1mm] u(t) \esmat\,dt \right]$.  For any control $u(t)$, $t \in [0,\infty)$, and for a given $T>0$, let $x_u(t)$ be the state reached at time $t=T$ starting from initial condition $x(0)$ and using the control $u(t)$, $t \in [0,T]$. Then 
 \begin{eqnarray*}
J_\infty^\ast & =& \!\! \min_{u(t),\, t \ge 0} \left[ \int_0^T \bmat{cc} \!\! x^\tra(t) \!\!  & u^\tra(t) \!\!  \emat  \!\! \bmat{cc}  \!\! Q \!\!  & S \!\!  \\ \!\!  S^\tra  \!\! & R \!\!  \emat \bmat{c}  \!\! x(t) \!\!  \\ \!\!  u(t) \!\!  \emat\,dt \right. \\
&& \qquad \left. + \int_T^\infty \bmat{cc} \!\!  x^\tra(t) \!\!  & u^\tra(t) \!\!  \emat \!\!  \bmat{cc} \!\!  Q  \!\! & S \!\!  \\  \!\!  S^\tra \!\!  & R  \!\! \emat \!\!  \bmat{c}  \!\! x(t) \!\!  \\  \!\! u(t) \!\!  \emat\,dt \right] \\
   &=& \!\! \min_{u(t),\, t \in [0,T)}  \underbrace{ \left[ \int_0^T  \!\!  \bmat{cc} \!\!   x^\tra(t) & u^\tra(t) \!\!   \emat \Pi \bmat{c}  \!\! x(t)  \!\! \\  \!\! u(t) \!\!  \emat\,dt +x_u^\tra(T)\,\bar{X}\,x_u(T)\right]}_{J_{T,\bar{X}}(x_0,u)}, 
   \end{eqnarray*}
 where the latter is due to the principle of optimality. 
Thus, $u(t)$, $t \in [0,\infty)$ minimises $J_{\infty}$ if and only if $u(t)$, $t \in [0,T]$ minimises $J_{T,\bar{X}}(x_0,u)$ and $u(t)$, $t \in [T,\infty)$ is such that 
 \[
 \int_T^\infty \bmat{cc} x^\tra(t) & u^\tra(t) \emat \bmat{cc} Q & S \\ S^\tra & R \emat \bmat{c} x(t) \\ u(t) \emat\,dt=x_u^\tra(T)\,\bar{X}\,x_u(T).
 \]

 The set of controls that minimise $J_\infty$ are those, and only those, that minimise $J_{T,\bar{X}}$. The optimal cost of the latter problem is independent of how big the value of $T$ is selected. This optimal cost is achieved by using the controls given by (\ref{optcontr}), where $P_T(t)$ is constant and equal to $\bar{X}$, since $\bar{X}$ is a stationary solution of (\ref{grde1}-\ref{number}) and $H=\bar{X}$.
 \endproof
 
 We conclude this section with a result that links the existence of solutions of the generalised Riccati equation with a geometric identity involving the smallest input containing subspace $\gS^\star$ and the largest reachability output-nulling subspace $\gR^\star$ of the underlying system, i.e., of the quadruple $(A,B,C,D)$, where $C$ and $D$ are matrices of suitable sizes such that
 \[
 \Pi=\bmat{c} C^\tra \\ D^\tra \emat \bmat{cc} C & D \emat.
 \]
 For more details of the underlying geometric concepts of input-containing and output-nulling subspaces, we refer to the monograph \cite{Trentelman-SH-01}. 
 
 {\begin{proposition}
Let CGCARE($\Sigma$) admit a solution $X=X^\tra$. Then, $\gS^\star=\gR^\star$.
\end{proposition}
}
\proof
Let  $X=X^\tra$ be a solution of CGCARE($\Sigma$). Observe also that CGCARE($\Sigma$) can be re-written as
\begin{equation}
\label{cgdare12}
\left\{ \begin{array}{ll} X\,A_0+A_0^\tra\,X-X\,B\,R^\dagger\,B^\tra\,X+Q_0=0 \\
\ker R \subseteq \ker X\,B \end{array} \right.
\end{equation}
where $A_0 \defi A-B\,R^\dagger S^\tra$ and $Q_0 \defi Q-S\,R^\dagger S^\tra$. Recall that $G=I_m-R^\dagger R$, so that $B_2 = B\,G$, and (\ref{cgdare12}) becomes
\begin{equation}
\label{cgdare13}
\left\{ \begin{array}{ll} X\,A_0+A_0^\tra\,X-X\,B\,R^\dagger\,B^\tra\,X+Q_0=0 \\
X\,B\,G=0 \end{array} \right.
\end{equation}
It is easy to see that $\ker X \subseteq \ker Q_0$. 
Indeed, by multiplying the first of (\ref{cgdare13}) on the left by $\xi^\tra$ and on the right by $\xi$, where $\xi \in \ker X$, we get $\xi^\tra\,Q_0\,\xi=0$. However, $Q_0$ is positive semidefinite, being the generalised Schur complement of $Q$ in $\Pi$. Hence, $Q_0\,\xi=0$, which implies $\ker X \subseteq \ker Q_0$. 
Since $X\,B\,G=0$, we get also $Q_0\,B\,G=0$. By post-multiplying the first of (\ref{cgdare13}) by a vector $\xi \in \ker X$ we find $X\,A_0\,\xi=0$, which says that $\ker X$ is $A_0$-invariant. This means that $\ker X$ is an $A_0$-invariant subspace containing the image of $B\,G$. Then, the reachable subspace of the pair $(A_0,BG)$, denoted by $\gR(A_0,B\,G)$, which is the smallest $A_0$-invariant subspace containing the image of $B\,G$, is contained in $\ker X$, i.e., $\gR(A_0,B\,G)\subseteq \ker X$. Therefore also $\gR(A_0,B\,G)\subseteq \ker Q_0$. Notice that $Q_0$ can be written as $C_0^\tra\,C_0$, where $C_0\defi C-D\,R^\dagger S^\tra$. Indeed,
\begin{eqnarray*}
C_0^\tra\,C_0 &=& C^\tra C-C^\tra D R^\dagger S^\tra-S R^\dagger D^\tra\,C+S R^\dagger D^\tra D R^\dagger S^\tra \\ 
&= &Q-S R^\dagger S-S R^\dagger S^\tra+S R^\dagger S^\tra=Q_0.
\end{eqnarray*}
Consider the two quadruples $(A,B,C,D)$ and $(A_0,B,C_0,D)$. We observe that the second is obtained directly from the first by applying the feedback input $u(t)=-R^\dagger S\,x(t)+v(t)$. We denote by $\gV^\star$, $\gR^\star$ the largest output-nulling and reachability subspace of $(A,B,C,D)$, and by $\gS^\star$ the smallest input-containing subspace of $(A,B,C,D)$. Likewise, we denote by $\gV_0^\star$, $\gR_0^\star$, $\gS_0^\star$ the same subspaces relative to the quadruple $(A_0,B,C_0,D)$. Thus, $\gV^\star= \gV^\star_0$, 
$\gR^\star= \gR^\star_0$, and $\gS^\star= \gS^\star_0$. 
 The first two identities are obvious, as output-nulling subspaces can be made invariant under state-feedback transformations and reachability is invariant under the same transformation. The third follows from \cite[Theorem 8.17]{Trentelman-SH-01}.
 There holds $\gR^\star=\gR(A_0,B\,G)$. Indeed, consider a state $x_1\in \gR(A_0,B\,G)$. There exists a control function $u$ driving the state from the origin to $x_1$, and we show that this control keeps the output at zero. Since $\ima(B\,G)=B\,\ker D$, such control can be chosen to satisfy $D\,u(t)=0$ for all $t\ge 0$. Moreover, as we have already seen, from $Q_0=C_0^\tra C_0$ and $\gR(A,BG)=\gR(A_0,BG)$ we have $C_0\,\gR(A_0,B\,G)=0$ since $\gR(A,B\,G)$ lies in $\ker Q_0$. Therefore, the output is identically zero. This implies that $\gR(A_0,B\,G) \subseteq \gR^\star$. However, the reachability subspace of $(A_0,B,C_0,D)$ cannot be greater than $\gR(A_0,B\,G)$, since 
 $D^\tra C_0=D^\tra (I_m-D\,(D^\tra D)^{\dagger}D^\tra)C=0$. Therefore, such control must necessarily render the output non-zero. The same argument can be used to prove that $\gS^\star=\gR(A_0,B\,G)$, where distributions can also be used in the allowed control, since $\gR(A,BG)$ represents also the set of states that are reachable from the origin using distributions in the control law \cite[p. 183]{Trentelman-SH-01}. Hence, $\gS^\star=\gR^\star$.
 \endproof

 \section{Concluding remarks}
 In this paper we established a new theory that showed that, when the CGCARE($\Sigma$) admits solutions, the corresponding singular LQ problem admits an impulse-free solution, and the optimal control can be expressed in terms of a state feedback.
 A very interesting question, which is currently being investigated by the authors, is the converse implication of this statement: when the singular LQ problem admits a regular solution for all initial states $x_0 \in \real^n$, does the CGCARE($\Sigma$) admit at least one symmetric positive semidefinite solution?
 At this stage we can only conjecture that this is the case, on the basis of some preliminary work, but the issue is indeed an open and interesting one.

\end{document}